\documentclass[12pt,a4paper]{article}

\usepackage{amsmath,amsthm,amssymb}
\usepackage{mathtools}
\usepackage{hyperref}
\usepackage{geometry}
\newtheorem{theorem}{Theorem}[section]
\newtheorem{proposition}[theorem]{Proposition}
\newtheorem{lemma}[theorem]{Lemma}

\newtheorem{remark}[theorem]{Remark}
\newtheorem{definition}[theorem]{Definition}

\newcommand{\R}{\mathbb{R}}
\newcommand{\E}{\mathbb{E}}

\begin{document}

\title{Stochastic differential delay equations with random delay\\
driven by fractional Brownian motion}

\author{Carles Rovira%
  \thanks{Facultat de Matem\`atiques i Inform\`atica,
    Universitat de Barcelona, Gran Via 585, 08007 Barcelona, Spain.
    E-mail: \texttt{Carles.Rovira@ub.edu}.
    Supported by grant PID2025-173091NB-I00
    funded by MCIN/AEI/10.13039/501100011033.}}

\maketitle



\begin{abstract}
We prove an existence and uniqueness result for stochastic differential
delay equations driven by a fractional Brownian motion with Hurst
parameter $H > \tfrac{1}{2}$, where the delay is a random process.
Extending the pathwise Riemann--Stieltjes approach of Ferrante and Rovira
\cite{FR2010}, we show that if the random delay $\tau(\cdot,\omega)$ has
almost surely $\gamma$-H\"older continuous trajectories with
$\gamma\delta > (1-H)/H$, where $\delta\in(0,1]$ is the H\"older exponent of
the gradient of the diffusion coefficient, then the equation has a unique
solution with paths in
$C^{1-\alpha}(-r_0,T;\R^d)\subset W_0^{\alpha,\infty}(-r_0,T;\R^d)$ for
appropriate $\alpha\in(1-H,\tfrac12)$.
We also prove that the solution has moments of all orders when
$\gamma>1/(2H)$ and the H\"older constant of the delay is bounded.
\end{abstract}

\section{Introduction}

Stochastic differential delay equations (SDDEs) arise naturally in
mathematical finance, population dynamics, control theory, and
neuroscience, wherever the present state of a system depends on its
past. 

In \cite{FR2006}, Ferrante and Rovira proved existence, uniqueness and
moment estimates for the equation
\begin{equation}\label{eq:fixed}
  X(t) = \eta(0) + \int_0^t b(s,X(s))\,ds
  + \int_0^t \sigma(s, X(s-r))\,dW_s^H,
  \quad t\in(0,T],
\end{equation}
with a fixed delay $r>0$, using the Riemann--Stieltjes integral
in the sense of Z\"ahle \cite{Zahle1998}.
In \cite{FR2010} they studied the convergence of solutions of
\eqref{eq:fixed} as $r\to 0$.

But it is well-known that sometimes  delays are not constants, for instance in neural networks \cite{Toub}, in population biology and
epidemiology \cite{Boukanjime2020,Martcheva2025}, in mathematical finance
\cite{CarteaSanchezBetancourt2023}, or in dynamical systems
\cite{GarridoAtienza2011}.

Our goal  is to extend the results of \cite{FR2010}
to the equation
\begin{equation}\label{eq:main}
\begin{aligned}
  X(t) &= \eta(0)
    + \int_0^t b(s,X(s), X(s-\tau(s)))\,ds\\
    &\quad+ \int_0^t \sigma(s, X(s), X(s-\tau(s)))\,dW_s^H,
    && t\in(0,T],\\
  X(t) &= \eta(t), && t\in[-r_0,0],
\end{aligned}
\end{equation}
where $\tau:[0,T]\times\Omega\to[0,r_0]$ is a random delay
with almost surely H\"older-continuous trajectories,
$r_0>0$ is the maximum delay, and $\eta:[-r_0,0]\to\R^d$ is the initial path.
Here both coefficients $b$ and $\sigma$ depend pointwise on the
current value $X(s)$ and the delayed value $X(s-\tau(s))$.

When the driving noise is a standard Brownian motion, the theory on delays, based
on It\^o calculus, is by now
well developed.  Stochastic functional differential equations with
general path-dependent coefficients are treated in the monographs of
Mohammed~\cite{Mohammed1984,Mohammed1998} and Mao~\cite{Mao1997}.
Random or state-dependent delays in the It\^o setting have also attracted
attention.  Random delays appear in Mohammed~\cite{Mohammed1984} and in
Kolmanovskii and Myshkis~\cite{KolMys}.  Mao, Shen and
Yuan~\cite{MaoShenYuan2008} study the almost sure asymptotic stability of
neutral stochastic differential delay equations with Markovian switching.

For the fractional Brownian case, the literature on delays is much
thinner.  Ferrante and Rovira~\cite{FR2006,FR2010} study delay equations by
classical methods, using a pathwise Riemann--Stieltjes approach.
Neuenkirch, Nourdin and Tindel~\cite{NeuenkirchNourdinTindel2008}
treat delay equations driven by rough paths, which cover fBm with
$H>\tfrac13$, using the Lyons rough-path theory \cite{Lyons1994}, and
Le\'on and Tindel~\cite{LeonTindel} study Malliavin calculus for fractional
delay equations with fixed delay.  As far as the author knows, it is the first study on random delays in the fBm setting.

In this paper, since $H>\tfrac{1}{2}$, the integral $\int_0^t \sigma\,dW^H$ is defined
pathwise as a Riemann--Stieltjes integral; see Young~\cite{Young1936},
Z\"ahle~\cite{Zahle1998} and Nualart--R\u{a}\c{s}canu~\cite{NR2002}.
For each fixed $\omega$, $W^H(\cdot,\omega)$ is a $(H-\varepsilon)$-H\"older
function, while $\tau(\cdot,\omega)$ is assumed to be $\gamma$-H\"older
with $\gamma\delta>(1-H)/H$.
Under these assumptions we work pathwise: for each $\omega$ we
solve a deterministic integral equation with the specific
realizations of $(W^H(\cdot,\omega),\tau(\cdot,\omega))$.
As in \cite{NR2002,FR2010}, the stochastic result is then obtained
by applying the deterministic one for almost every $\omega$.

The main difference with respect
to \cite{NR2002,FR2010} is  that we work with H\"older seminorms, following
\cite{HuNualart2007}, and we use a contraction argument only on short
time intervals.
In the above papers the solution is obtained as the fixed point of a
contraction on the whole interval $[0,T]$, with respect to the weighted
norms of the spaces $W_0^{\alpha,\infty}$, but
when the delay is variable, we do not know how to control the delayed term
$X(s-\tau(s))$ in these norms.

The paper is organized as follows.
Section~\ref{sec:prelim} introduces the function spaces, defines the
Riemann--Stieltjes integral, and states the hypotheses on the
coefficients.
Section~\ref{sec:main} states the main result, Theorem~\ref{thm:main}.
Section~\ref{sec:composition} contains the key new Lemma~\ref{lem:comp},
which controls the regularity of the composition $s\mapsto f(s-\tau(s))$.
Section~\ref{sec:det} proves existence and uniqueness for the deterministic
version of~\eqref{eq:main}, together with an explicit a priori bound.
Section~\ref{sec:stoch} applies these results to the stochastic equation
and completes the proof of Theorem~\ref{thm:main}.

\section{Preliminaries}\label{sec:prelim}

\subsection{Function spaces}\label{ss:spaces}

Fix parameters $r_0>0$, $T>0$ and $\alpha\in(0,\tfrac{1}{2})$.
We first recall the H\"older spaces.

\begin{definition}
For $\lambda\in(0,1]$ and an interval $I\subset\R$, the H\"older space
$C^\lambda(I;\R^d)$ consists of all continuous functions $f:I\to\R^d$
such that the seminorm
\[
  [f]_{C^\lambda(I)} := \sup_{\substack{s,t\in I\\s\neq t}}
  \frac{|f(t)-f(s)|}{|t-s|^\lambda} < \infty.
\]
We write $C^\lambda$ for $C^\lambda([-r_0,T];\R^d)$ when the interval and
target are clear, and $\|f\|_{C^\lambda(I)} := \|f\|_{L^\infty(I)}
+ [f]_{C^\lambda(I)}$.
\end{definition}

We shall also need the spaces $W_0^{\alpha,\infty}$ of \cite{NR2002}.

\begin{definition}\label{def:Walpha}
For $f:[-r_0,T]\to\R^d$ measurable, define the
norm
\begin{equation}\label{eq:weightnorm}
  \|f\|_{\alpha,r_0}
  := \sup_{t\in[-r_0,T]}
  \left(|f(t)| + \int_{-r_0}^{t}
  \frac{|f(t)-f(s)|}{(t-s)^{\alpha+1}}\,ds\right).
\end{equation}
We denote
\begin{equation*}
  W_0^{\alpha,\infty}(-r_0,T;\R^d)
  := \bigl\{f:[-r_0,T]\to\R^d \text{ measurable}:
  \|f\|_{\alpha,r_0}<\infty\bigr\}.
\end{equation*}
For functions defined only on $[0,T]$, we write
$W_0^{\alpha,\infty}(0,T;\R^d)$ with the same formula but integration
from $0$, and we denote the corresponding norm by $\|f\|_{\alpha,0}$.
\end{definition}

\begin{remark}\label{rem:singint}
Notice that if
$f$ is $\beta$-H\"older continuous with $\beta>\alpha$, then
\[
  \int_{-r_0}^t \frac{|f(t)-f(s)|}{(t-s)^{\alpha+1}}\,ds
  \leq [f]_{C^\beta}\int_{-r_0}^t (t-s)^{\beta-\alpha-1}\,ds
  = \frac{[f]_{C^\beta}\,(t+r_0)^{\beta-\alpha}}{\beta-\alpha} < \infty.
\]
In particular, $C^{1-\alpha}([-r_0,T];\R^d)\subset W_0^{\alpha,\infty}(-r_0,T;\R^d)$
with
\begin{equation}\label{eq:holder-embed}
  \|f\|_{\alpha,r_0} \leq C_{\alpha,T,r_0}\,\|f\|_{C^{1-\alpha}([-r_0,T])}.
\end{equation}
\end{remark}

\subsection{The Riemann--Stieltjes integral}\label{ss:RS}

The Riemann--Stieltjes integral $\int_0^t f(s)\,dg_s$ is not well-defined
for all pairs $(f,g)$ of H\"older functions. The theory of Nualart and
R\u{a}\c{s}canu~\cite{NR2002}, following Z\"ahle~\cite{Zahle1998},
provides sufficient conditions via fractional calculus.

For the definition of the integral $\int f\,dg$ and of the constant
$\Lambda_\alpha(g)$ associated with a function $g:[0,T]\to\R^m$ and
$\alpha\in(0,\tfrac12)$ we refer to \cite{NR2002}; here we only recall the
two facts that will be used.

\begin{remark}\label{rem:Lambda-holder}
If $g\in C^{\beta'}([0,T];\R^m)$ with $\beta'>1-\alpha$, then, as shown
in \cite{NR2002}, $\Lambda_\alpha(g)<\infty$ and, more precisely,
\begin{equation}\label{eq:Lambda-holder}
  \Lambda_\alpha(g)\leq C_{\alpha,\beta',T}\,[g]_{C^{\beta'}([0,T])}.
\end{equation}
In particular, for $g=W^H$ with $H>\tfrac12$, choosing
$\alpha\in(1-H,\tfrac12)$ ensures $\Lambda_\alpha(W^H)<\infty$ a.s.
\end{remark}

The estimate that we shall use is the following one, from
\cite[p.~61]{NR2002}.
Let $g$ satisfy $\Lambda_\alpha(g)<\infty$, let $0\leq u<s\leq T$ and let
$f:[u,s]\to\R^{d\times m}$ be measurable with
\[
  \int_u^s\bigl(|f(r)|(r-u)^{-\alpha}
  +\int_u^r|f(r)-f(w)|(r-w)^{-\alpha-1}\,dw\bigr)dr<\infty.
\]
Then the integral $\int_u^s f\,dg$ is well defined and
\begin{equation}\label{eq:est4}
  \left|\int_u^s f(r)\,dg_r\right| \leq \Lambda_\alpha(g)
    \int_u^s\left( \frac{|f(r)|}{(r-u)^\alpha}
    + \alpha\int_u^r \frac{|f(r)-f(w)|}{(r-w)^{\alpha+1}}\,dw\right)dr .
\end{equation}
Moreover the integral is additive with respect to the interval of
integration: $\int_u^s f\,dg=\int_u^v f\,dg+\int_v^s f\,dg$ for $u<v<s$.
When $g=(g^1,\dots,g^m)$ takes values in $\R^m$ and $f$ in
$\R^{d\times m}$, the integral is understood as
$\int f\,dg=\sum_{j=1}^m\int f^{\cdot j}\,dg^j$, where $f^{\cdot j}$ is the
$j$-th column of $f$, and we set
$\Lambda_\alpha(g):=\sum_{j=1}^m\Lambda_\alpha(g^j)$.
Since $|f^{\cdot j}|\leq|f|$, the estimate~\eqref{eq:est4} then holds as
stated, by applying the scalar estimate to each term of the sum.

\subsection{Hypotheses}\label{ss:hyp}

For $(x,y)\in\R^d\times\R^d$ we use the norm $|(x,y)|:=|x|+|y|$, we write
$\nabla\sigma:=\nabla_{(x,y)}\sigma$ for the gradient of $\sigma$ with
respect to both spatial arguments, and $\|\cdot\|$ for the
corresponding operator norm.

\begin{itemize}
\item[\textbf{(H1)}]
$\sigma:[0,T]\times\R^d\times\R^d\to\R^{d\times m}$ is measurable,
$\sigma(t,\cdot,\cdot)$ is $C^1$ for each $t$,
and there exist constants $\beta,\delta\in(0,1]$, $M_0>0$ and, for each
$N\geq 0$, a constant $M_N>0$ such that for all
$s,t\in[0,T]$:
\begin{enumerate}
  \item[(i)] $|\sigma(t,x,y)-\sigma(t,x',y')|\leq M_0(|x-x'|+|y-y'|)$
    for all $x,y,x',y'\in\R^d$;
  \item[(ii)] $\|\nabla\sigma(t,x,y)-\nabla\sigma(t,x',y')\|
    \leq M_N\bigl(|x-x'|+|y-y'|\bigr)^\delta$
    for all $x,y,x',y'$ with $|x|,|y|,|x'|,|y'|\leq N$;
  \item[(iii)] $|\sigma(t,x,y)-\sigma(s,x,y)|
    + \|\nabla\sigma(t,x,y)-\nabla\sigma(s,x,y)\|
    \leq M_0|t-s|^\beta$ for all $x,y\in\R^d$.
\end{enumerate}

\item[\textbf{(H2)}]
$b:[0,T]\times\R^d\times\R^d\to\R^d$ is measurable and
there exist a function $b_0\in L^{1/\alpha}(0,T)$, a constant $L_0\geq 0$ and,
for each $N\geq 0$, a constant $L_N>0$ such that for all
$t\in[0,T]$:
\begin{enumerate}
  \item[(i)] $|b(t,x,y)-b(t,x',y')|\leq L_N\bigl(|x-x'|+|y-y'|\bigr)$
  for all $x,y,x',y'\in\R^d$ such that $|x|,|y|,|x'|,|y'|\leq N$;
  \item[(ii)] $|b(t,x,y)|\leq L_0\bigl(|x|+|y|\bigr)+b_0(t)$
  for all $x,y\in\R^d$.
\end{enumerate}

\item[\textbf{(H$\tau$)}]
The random delay $\tau:[0,T]\times\Omega\to[0,r_0]$
is a measurable stochastic process, and there exist
$\gamma\in(0,1]$ with
\[
  \gamma\,\delta > \frac{1-H}{H}
\]
and a random variable $C_\tau\geq 0$, finite almost surely,
such that, almost surely,
\begin{equation}\label{eq:tau-holder}
  |\tau(t,\omega)-\tau(s,\omega)| \leq C_\tau(\omega)\,|t-s|^\gamma,
  \quad\forall\, s,t\in[0,T].
\end{equation}
\end{itemize}

\begin{remark}\label{rem:H1}
As in \cite{NR2002,FR2010}, the Lipschitz condition \textbf{(H1)(i)} is
global.
It implies the linear growth bound
\begin{equation}\label{eq:lin-growth}
  |\sigma(t,x,y)|\leq K_0+M_0(|x|+|y|),\qquad
  K_0:=|\sigma(0,0,0)|+M_0T^\beta ,
\end{equation}
so no separate growth hypothesis on $\sigma$ is needed.
\end{remark}

\begin{remark}\label{rem:tau-examples}
When $\delta=1$, for instance if $\sigma(t,\cdot,\cdot)$ is $C^2$ with
locally bounded second derivatives, the condition in \textbf{(H$\tau$)}
reads $\gamma>(1-H)/H$.
It is satisfied, for example, by any Lipschitz delay, for which $\gamma=1$,
or by $\tau(t)=\psi(Y_t)$ with $\psi:\R\to[0,r_0]$ Lipschitz and $Y$ a
process with $\gamma$-H\"older paths.
\end{remark}

\section{Main result}\label{sec:main}

Let $(\Omega,\mathcal{F},\mathbb{P})$ be a complete probability space and
$W^H=(W^{H,j})_{j=1}^m$ independent fractional Brownian motions with
Hurst parameter $H>\tfrac{1}{2}$.
Recall that $W^{H,j}$ is a centered Gaussian process with covariance
$\E[W^{H,j}(t)W^{H,j}(s)]=\tfrac{1}{2}(t^{2H}+s^{2H}-|t-s|^{2H})$.
We consider equation~\eqref{eq:main} under the hypotheses of
Section~\ref{ss:hyp}. The main result of the paper is the following.

\begin{theorem}\label{thm:main}
Let $H>\tfrac12$. Assume \textbf{(H1)} with $\beta>1-H$, \textbf{(H2)}
and \textbf{(H$\tau$)}, and set
\[
  \alpha_0:=\min\Bigl\{\tfrac12,\ \beta,\ \tfrac{\delta\gamma}{1+\delta\gamma}\Bigr\}
  >1-H .
\]
Fix $\alpha\in(1-H,\alpha_0)$, assume $b_0\in L^{1/\alpha}(0,T)$, and let
$\eta$ be a stochastic process with trajectories in
$C^{1-\alpha}([-r_0,0];\R^d)$ almost surely.
\begin{enumerate}
  \item[\emph{(a)}]
  There exists a unique stochastic process $X$ with trajectories in
  $C^{1-\alpha}([-r_0,T];\R^d)$ almost surely which is a solution of the
  stochastic equation~\eqref{eq:main}. In particular,
  \[
    X\in L^0\bigl(\Omega,\mathcal F,\mathbb P;\,W_0^{\alpha,\infty}(-r_0,T;\R^d)\bigr).
  \]
  \item[\emph{(b)}]
  Assume in addition that
  \begin{equation}\label{eq:moment-cond}
    \gamma>\frac{1}{2H},\qquad \alpha<1-\frac{1}{2\gamma},\qquad
    C_\tau\in L^\infty(\Omega),\qquad
    \E\bigl[\|\eta\|_{C^{1-\alpha}([-r_0,0])}^p\bigr]<\infty\ \ \forall p\geq1 .
  \end{equation}
  Then, for all $p\geq 1$,
  \[
    \E\bigl[\|X\|_{C^{1-\alpha}([-r_0,T])}^p\bigr]<\infty,
    \qquad\text{and hence}\qquad
    \E\bigl[\|X\|_{\alpha,r_0}^p\bigr]<\infty .
  \]
\end{enumerate}
\end{theorem}

\begin{remark}
The inequality $\alpha_0>1-H$ holds because $\beta>1-H$ and
because $\delta\gamma/(1+\delta\gamma)>1-H$ is equivalent to
$\delta\gamma>(1-H)/H$, which is assumed in \textbf{(H$\tau$)}.
For a Lipschitz delay, that is, for $\gamma=1$, the first two conditions
in~\eqref{eq:moment-cond} hold automatically.
\end{remark}

The proof of Theorem~\ref{thm:main} is given in Section~\ref{sec:stoch}.
It relies on the composition lemma of Section~\ref{sec:composition} and on
the deterministic results of Section~\ref{sec:det}.

\section{The composition lemma}\label{sec:composition}

The following lemma is the key new ingredient of this paper.

\begin{lemma}\label{lem:comp}
Let $\alpha\in(0,\tfrac{1}{2})$, $\gamma>\alpha/(1-\alpha)$,
$r_0,T>0$.
Suppose $f\in C^{1-\alpha}(-r_0,T;\R^d)$
and $\tau:[0,T]\to[0,r_0]$ is $\gamma$-H\"older with constant
$[\tau]_\gamma := [{\tau}]_{C^\gamma([0,T])}$.
Define $f^\tau:[0,T]\to\R^d$ by $f^\tau(s):= f(s-\tau(s))$.
Then $f^\tau\in W_0^{\alpha,\infty}(0,T;\R^d)$; more precisely:
\begin{enumerate}
  \item[\emph{(a)}]
  For all $s,t\in[0,T]$,
  \begin{equation}\label{eq:comp-holder}
    |f^\tau(t)-f^\tau(s)|
    \leq [f]_{C^{1-\alpha}}\,(T^{1-\gamma}+[\tau]_\gamma)^{1-\alpha}\,
    |t-s|^{\gamma(1-\alpha)}.
  \end{equation}
  \item[\emph{(b)}]
  For all $t\in[0,T]$,
  \begin{equation}\label{eq:comp-bound}
    \int_0^t\frac{|f^\tau(t)-f^\tau(s)|}{(t-s)^{\alpha+1}}\,ds
    \leq \frac{(T^{1-\gamma}+[\tau]_\gamma)^{1-\alpha}\,[f]_{C^{1-\alpha}}\,
    T^{\gamma(1-\alpha)-\alpha}}{\gamma(1-\alpha)-\alpha}.
  \end{equation}
  \item[\emph{(c)}]
  $\|f^\tau\|_{\alpha,0}
  \leq \|f\|_{\alpha,r_0}
  + \dfrac{(T^{1-\gamma}+[\tau]_\gamma)^{1-\alpha}\,[f]_{C^{1-\alpha}}\,
  T^{\gamma(1-\alpha)-\alpha}}{\gamma(1-\alpha)-\alpha}$.
\end{enumerate}
\end{lemma}

\begin{proof}
We first obtain the uniform bound.
Since $t-\tau(t)\in[-r_0,T]$, we have $|f^\tau(t)|\leq \|f\|_{L^\infty(-r_0,T)}
\leq \|f\|_{\alpha,r_0}$.

We now prove the H\"older estimate~(a).
Since $f\in C^{1-\alpha}(-r_0,T)$,
\begin{equation}\label{eq:f-holder-comp}
  |f^\tau(t)-f^\tau(s)|
  = |f(t-\tau(t))-f(s-\tau(s))|
  \leq [f]_{C^{1-\alpha}}\,\bigl|(t-\tau(t))-(s-\tau(s))\bigr|^{1-\alpha}.
\end{equation}
Since $\gamma\leq 1$ and $|t-s|\leq T$, we have
$|t-s|=|t-s|^{1-\gamma}\cdot|t-s|^\gamma \leq T^{1-\gamma}|t-s|^\gamma$,
so
\begin{align}
  \bigl|(t-\tau(t))-(s-\tau(s))\bigr|
  &\leq |t-s|+|\tau(t)-\tau(s)| \notag\\
  &\leq T^{1-\gamma}|t-s|^\gamma+[\tau]_\gamma|t-s|^\gamma
   = (T^{1-\gamma}+[\tau]_\gamma)|t-s|^\gamma. \label{eq:tau-shift}
\end{align}
Inserting into~\eqref{eq:f-holder-comp} gives~(a).

To obtain~(b), we use part~(a):
\[
  \int_0^t \frac{|f^\tau(t)-f^\tau(s)|}{(t-s)^{\alpha+1}}\,ds
  \leq [f]_{C^{1-\alpha}}\,(T^{1-\gamma}+[\tau]_\gamma)^{1-\alpha}
  \int_0^t (t-s)^{\gamma(1-\alpha)-\alpha-1}\,ds.
\]
The exponent satisfies $\gamma(1-\alpha)-\alpha-1>-1$ if and only if
$\gamma(1-\alpha)>\alpha$, i.e.\ $\gamma>\alpha/(1-\alpha)$, which holds by
assumption.  The integral equals
$t^{\gamma(1-\alpha)-\alpha}/(\gamma(1-\alpha)-\alpha)
\leq T^{\gamma(1-\alpha)-\alpha}/(\gamma(1-\alpha)-\alpha)$,
giving~(b).

Finally, combining the uniform bound with~(b) we obtain~(c).
\end{proof}

\section{The deterministic equation}\label{sec:det}

In this section we fix $\omega\in\Omega$ and work deterministically
with the realizations $g:=W^H(\cdot,\omega)$ and $\tau:=\tau(\cdot,\omega)$.
Thus $g:[0,T]\to\R^m$ is a fixed function with $\Lambda_\alpha(g)<\infty$
and $\tau:[0,T]\to[0,r_0]$ is a fixed $\gamma$-H\"older function.

\subsection{Notation}\label{ss:notation}

Throughout this section \textbf{(H1)} and \textbf{(H2)} hold, and
\begin{equation}\label{eq:alpha0}
  \alpha\in(0,\alpha_0),\qquad
  \alpha_0 := \min\Bigl\{\tfrac{1}{2},\,\beta,\,
  \tfrac{\delta\gamma}{1+\delta\gamma}\Bigr\}.
\end{equation}
We write
\begin{gather*}
  \mu:=1-\alpha,\qquad
  \Lambda:=\Lambda_\alpha(g),\qquad
  \mathcal{T}:=T^{1-\gamma}+[\tau]_\gamma,\\
  \kappa:=\min\{\beta,\ \delta\gamma\mu\},\qquad
  c_0:=\frac{1}{\kappa-\alpha},\qquad
  C_1:=\|b_0\|_{L^{1/\alpha}(0,T)}+\frac{K_0+M_0}{16M_0},
\end{gather*}
where $K_0$ is the constant in~\eqref{eq:lin-growth}.
By~\eqref{eq:alpha0} we have $\beta>\alpha$ and $\delta\gamma\mu>\alpha$,
so $\kappa>\alpha$ and $c_0<\infty$; note also that every exponent
$\theta\in\{1,\beta,\mu,\gamma\mu,\delta\gamma\mu\}$ satisfies
$1/(\theta-\alpha)\leq c_0$.
For an interval $I$ we write $\|f\|_{\infty,I}:=\sup_I|f|$ and
$[f]_{\mu,I}:=[f]_{C^\mu(I)}$.
For $0<b\leq T$ and $f:[-r_0,b]\to\R^d$ we set
\[
  f^\tau(s):=f(s-\tau(s)),\qquad \Phi_f(s):=\bigl(f(s),f^\tau(s)\bigr),
  \qquad s\in[0,b].
\]
By Lemma~\ref{lem:comp}(a) and
$|s-u|^\mu\leq T^{(1-\gamma)\mu}|s-u|^{\gamma\mu}\leq\mathcal{T}^\mu|s-u|^{\gamma\mu}$,
for all $u,s\in[0,b]$:
\begin{equation}\label{eq:Phi-holder}
  |f^\tau(s)-f^\tau(u)|\leq \mathcal{T}^\mu[f]_{\mu,[-r_0,b]}|s-u|^{\gamma\mu},
  \qquad
  |\Phi_f(s)-\Phi_f(u)|\leq 2\mathcal{T}^\mu[f]_{\mu,[-r_0,b]}|s-u|^{\gamma\mu}.
\end{equation}
Given $0\leq a<b\leq T$ and $f\in C^\mu([-r_0,b];\R^d)$ we define, for
$t\in[a,b]$,
\begin{align}
  F_a(f)(t) &:= \int_a^t b(r, f(r), f^\tau(r))\,dr,
  \label{eq:Fb}\\
  G_a(f)(t) &:= \int_a^t \sigma(r, f(r), f^\tau(r))\,dg_r.
  \label{eq:Gsigma}
\end{align}

\subsection{Estimates on the intervals}

We start with two basic lemmas.
The first lemma is the form of the estimate~\eqref{eq:est4} that we shall
use, for H\"older continuous integrands.

\begin{lemma}\label{lem:young}
Let $0\leq u<s\leq T$ and let $\varphi:[u,s]\to\R^{d\times m}$ be bounded and
such that
\[
  |\varphi(r)-\varphi(w)|\leq\sum_{i=1}^k a_i\,|r-w|^{\theta_i},
  \qquad r,w\in[u,s],
\]
for some $a_i\geq0$ and $\theta_i\in(\alpha,1]$.
Then $\int_u^s\varphi\,dg$ is well defined and
\begin{equation}\label{eq:young}
  \left|\int_u^s\varphi(r)\,dg_r\right|
  \leq \Lambda\,(s-u)^{1-\alpha}\left(
  \frac{\|\varphi\|_{\infty,[u,s]}}{1-\alpha}
  +\sum_{i=1}^k\frac{a_i\,(s-u)^{\theta_i}}{\theta_i-\alpha}\right).
\end{equation}
\end{lemma}

\begin{proof}
We insert the assumptions into~\eqref{eq:est4}.
First, $$\int_u^s|\varphi(r)|(r-u)^{-\alpha}dr
\leq\|\varphi\|_{\infty,[u,s]}(s-u)^{1-\alpha}/(1-\alpha).$$
Second, for each $i$,
\[
  \alpha\int_u^s\!\!\int_u^r (r-w)^{\theta_i-\alpha-1}\,dw\,dr
  =\frac{\alpha\,(s-u)^{1+\theta_i-\alpha}}{(\theta_i-\alpha)(1+\theta_i-\alpha)}
  \leq\frac{(s-u)^{1-\alpha}(s-u)^{\theta_i}}{\theta_i-\alpha},
\]
because $\alpha<1\leq1+\theta_i-\alpha$.
\end{proof}

The second lemma collects the estimates for the increments of $\sigma$.

\begin{lemma}\label{lem:sigma}
Assume \textbf{(H1)} and write $z=(x,y)$, $|z|=|x|+|y|$.
For all $s,t\in[0,T]$:
\begin{enumerate}
  \item[\emph{(a)}] $|\sigma(t,z)-\sigma(s,z')|\leq M_0|t-s|^\beta+M_0|z-z'|$
  for all $z,z'$;
  \item[\emph{(b)}] $\bigl|[\sigma(t,z_1)-\sigma(s,z_1)]-[\sigma(t,z_3)-\sigma(s,z_3)]\bigr|
  \leq M_0|t-s|^\beta\,|z_1-z_3|$ for all $z_1,z_3$;
  \item[\emph{(c)}] if $z_i=(x_i,y_i)$ with $|x_i|,|y_i|\leq N$ for
  $i=1,\dots,4$, then
  \begin{equation}\label{eq:second-diff}
  \begin{aligned}
    &\bigl|\sigma(t,z_1)-\sigma(t,z_2)-\sigma(t,z_3)+\sigma(t,z_4)\bigr|\\
    &\qquad\leq M_0\,|z_1-z_2-z_3+z_4|
    +M_N\,|z_1-z_3|\,\bigl(|z_1-z_2|^\delta+|z_3-z_4|^\delta\bigr).
  \end{aligned}
  \end{equation}
\end{enumerate}
\end{lemma}

\begin{proof}
(a) follows from \textbf{(H1)(i)} and \textbf{(H1)(iii)}.
Note that \textbf{(H1)(i)} gives $\|\nabla\sigma(t,z)\|\leq M_0$ for all $z$.

(b) Apply the mean value theorem in integral form to the $C^1$ map
$\psi(z):=\sigma(t,z)-\sigma(s,z)$ on the segment $[z_3,z_1]$:
$$\psi(z_1)-\psi(z_3)=\int_0^1\nabla\psi(z_3+\theta(z_1-z_3))\,d\theta\,(z_1-z_3),$$
and by \textbf{(H1)(iii)} the gradient of $\psi$ is bounded in norm by
$M_0|t-s|^\beta$.
Here the gradient is taken with respect to both arguments $x$ and $y$,
which is why \textbf{(H1)(iii)} is stated for the full gradient.

(c) By the mean value theorem in integral form, applied in the variable
$z=(x,y)$,
\begin{align*}
  \sigma(t,z_1)-\sigma(t,z_3)
  &=\int_0^1\nabla\sigma\bigl(t,z_3+\theta(z_1-z_3)\bigr)\,d\theta\;(z_1-z_3),\\
  \sigma(t,z_2)-\sigma(t,z_4)
  &=\int_0^1\nabla\sigma\bigl(t,z_4+\theta(z_2-z_4)\bigr)\,d\theta\;(z_2-z_4).
\end{align*}
Subtracting, and adding and subtracting
$\int_0^1\nabla\sigma(t,z_4+\theta(z_2-z_4))\,d\theta\,(z_1-z_3)$, we get
\begin{align*}
  &\sigma(t,z_1)-\sigma(t,z_2)-\sigma(t,z_3)+\sigma(t,z_4)\\
  &\quad=\int_0^1\nabla\sigma\bigl(t,z_4+\theta(z_2-z_4)\bigr)\,d\theta\;
  \bigl[(z_1-z_3)-(z_2-z_4)\bigr]\\
  &\qquad+\int_0^1\Bigl[\nabla\sigma\bigl(t,z_3+\theta(z_1-z_3)\bigr)
  -\nabla\sigma\bigl(t,z_4+\theta(z_2-z_4)\bigr)\Bigr]d\theta\;(z_1-z_3).
\end{align*}
The first term is bounded by $M_0|z_1-z_2-z_3+z_4|$.
In the second term all the points involved lie in the convex set
$\{|x|\leq N,\ |y|\leq N\}$ and
\[
  \bigl|z_3+\theta(z_1-z_3)-z_4-\theta(z_2-z_4)\bigr|
  =\bigl|(1-\theta)(z_3-z_4)+\theta(z_1-z_2)\bigr|
  \leq|z_3-z_4|+|z_1-z_2| ,
\]
so by \textbf{(H1)(ii)} and $(p+q)^\delta\leq p^\delta+q^\delta$ the second
term is bounded by $$M_N|z_1-z_3|(|z_1-z_2|^\delta+|z_3-z_4|^\delta).$$
\end{proof}

We can now estimate the H\"older seminorms of $F_a(f)$ and $G_a(f)$ on a
short interval.

\begin{proposition}\label{prop:growth}
Let $0\leq a<b\leq T$ with $\ell:=b-a\leq1$ and
$f\in C^\mu([-r_0,b];\R^d)$.
Then $F_a(f),G_a(f)\in C^\mu([a,b];\R^d)$ and
\begin{align}
  [F_a(f)]_{\mu,[a,b]}
  &\leq 2L_0\,\ell^{\alpha}\,\|f\|_{\infty,[-r_0,b]}+\|b_0\|_{L^{1/\alpha}(0,T)},
  \label{eq:F-growth}\\
  [G_a(f)]_{\mu,[a,b]}
  &\leq c_0\Lambda\Bigl(K_0+M_0+2M_0\|f\|_{\infty,[-r_0,b]}
  +M_0\bigl(1+\mathcal{T}^\mu\bigr)\ell^{\gamma\mu}\,[f]_{\mu,[-r_0,b]}\Bigr).
  \label{eq:G-growth}
\end{align}
\end{proposition}

\begin{proof}
Let $a\leq u<s\leq b$.
By \textbf{(H2)(ii)} and H\"older's inequality,
\[
  |F_a(f)(s)-F_a(f)(u)|
  \leq 2L_0\|f\|_{\infty,[-r_0,b]}(s-u)+\|b_0\|_{L^{1/\alpha}}(s-u)^{1-\alpha},
\]
and $(s-u)\leq\ell^\alpha(s-u)^{\mu}$; this gives~\eqref{eq:F-growth}.

Put $\varphi(r):=\sigma(r,\Phi_f(r))$.
By~\eqref{eq:lin-growth}, $\|\varphi\|_{\infty,[a,b]}\leq K_0+2M_0\|f\|_{\infty,[-r_0,b]}$.
By Lemma~\ref{lem:sigma}(a) and~\eqref{eq:Phi-holder}, for $r,w\in[a,b]$,
\[
  |\varphi(r)-\varphi(w)|
  \leq M_0|r-w|^\beta+M_0[f]_{\mu,[-r_0,b]}|r-w|^\mu
  +M_0\mathcal{T}^\mu[f]_{\mu,[-r_0,b]}|r-w|^{\gamma\mu}.
\]
Lemma~\ref{lem:young} with exponents $\beta,\mu,\gamma\mu$, all larger than
$\alpha$, together with $s-u\leq\ell\leq1$ and the definition of $c_0$,
yields
\[
  \left|\int_u^s\varphi\,dg\right|
  \leq c_0\Lambda(s-u)^{\mu}\Bigl(K_0+2M_0\|f\|_{\infty,[-r_0,b]}+M_0
  +M_0(1+\mathcal{T}^\mu)\ell^{\gamma\mu}[f]_{\mu,[-r_0,b]}\Bigr),
\]
which is~\eqref{eq:G-growth}.
\end{proof}

Next we prove the corresponding Lipschitz estimate, for two functions that
coincide up to time $a$.

\begin{proposition}\label{prop:lip}
Let $0\leq a<b\leq T$ with $\ell:=b-a\leq1$, and let
$f,h\in C^\mu([-r_0,b];\R^d)$ be such that
\[
  f=h \text{ on } [-r_0,a],\qquad
  \|f\|_{\infty,[-r_0,b]},\|h\|_{\infty,[-r_0,b]}\leq N,\qquad
  [f]_{\mu,[-r_0,b]},[h]_{\mu,[-r_0,b]}\leq R .
\]
Then
\begin{align}
  [F_a(f)-F_a(h)]_{\mu,[a,b]}
  &\leq 2L_N\,\ell\,[f-h]_{\mu,[a,b]},\label{eq:F-lip}\\
  [G_a(f)-G_a(h)]_{\mu,[a,b]}
  &\leq c_0\Lambda\,\ell^{\gamma\mu}\Bigl(5M_0+M_0\mathcal{T}^\mu
  +8M_N\mathcal{T}^{\delta\mu}R^\delta\Bigr)[f-h]_{\mu,[a,b]}.\label{eq:G-lip}
\end{align}
\end{proposition}

\begin{proof}
Write $\Delta:=f-h$ and $d:=[\Delta]_{\mu,[a,b]}$.
Since $\Delta=0$ on $[-r_0,a]$,
\begin{equation}\label{eq:Delta-basic}
  \|\Delta\|_{\infty,[-r_0,b]}\leq \ell^\mu d,
  \qquad
  [\Delta]_{\mu,[-r_0,b]}=d .
\end{equation}
For the second identity, note that if $u\leq a<s$ then
$$|\Delta(s)-\Delta(u)|=|\Delta(s)-\Delta(a)|\leq d(s-a)^\mu\leq d(s-u)^\mu.$$
Consequently, by Lemma~\ref{lem:comp}(a) applied to $\Delta$,
for all $u,s\in[a,b]$,
\begin{equation}\label{eq:Delta-tau}
  |\Delta^\tau(s)|\leq\ell^\mu d,\qquad
  |\Delta^\tau(s)-\Delta^\tau(u)|\leq \mathcal{T}^\mu d\,|s-u|^{\gamma\mu}.
\end{equation}

By \textbf{(H2)(i)}, for $a\leq u<s\leq b$,
\[
  \bigl|(F_a(f)-F_a(h))(s)-(F_a(f)-F_a(h))(u)\bigr|
  \leq L_N\int_u^s\bigl(|\Delta(r)|+|\Delta^\tau(r)|\bigr)dr
  \leq 2L_N\ell^\mu d\,(s-u),
\]
and $(s-u)\leq\ell^\alpha(s-u)^\mu$; this gives~\eqref{eq:F-lip}.

Put $v(s):=\sigma(s,\Phi_f(s))-\sigma(s,\Phi_h(s))$, so that
$G_a(f)-G_a(h)=\int_a^\cdot v\,dg$.
By \textbf{(H1)(i)}, \eqref{eq:Delta-basic} and~\eqref{eq:Delta-tau},
\begin{equation}\label{eq:v-pointwise}
  |v(s)|\leq M_0\bigl(|\Delta(s)|+|\Delta^\tau(s)|\bigr)\leq 2M_0\,\ell^\mu d .
\end{equation}
For $u,s\in[a,b]$ we add and subtract
$\sigma(u,\Phi_f(s))-\sigma(u,\Phi_h(s))$ and obtain the exact identity
\begin{equation}\label{eq:v-decomp}
  v(s)-v(u)=J_1(s,u)+J_2(s,u),
\end{equation}
where
\begin{align}
  J_1(s,u)&:=
  \bigl[\sigma(s,\Phi_f(s))-\sigma(u,\Phi_f(s))\bigr]
  -\bigl[\sigma(s,\Phi_h(s))-\sigma(u,\Phi_h(s))\bigr],
  \label{eq:J1}\\
  J_2(s,u)&:=
  \sigma(u,\Phi_f(s))-\sigma(u,\Phi_f(u))
  -\sigma(u,\Phi_h(s))+\sigma(u,\Phi_h(u)).
  \label{eq:J2}
\end{align}

\smallskip
\noindent We can get a bound on $J_1$.
By Lemma~\ref{lem:sigma}(b) with $z_1=\Phi_f(s)$, $z_3=\Phi_h(s)$,
\begin{equation}\label{eq:J1-bound}
  |J_1(s,u)|\leq M_0|s-u|^\beta\bigl(|\Delta(s)|+|\Delta^\tau(s)|\bigr)
  \leq 2M_0\,\ell^\mu d\,|s-u|^\beta .
\end{equation}

\noindent On the other hand, to bound $J_2$, 
we apply Lemma~\ref{lem:sigma}(c) at time $u$ with
\[
  z_1=\Phi_f(s),\quad z_2=\Phi_f(u),\quad z_3=\Phi_h(s),\quad z_4=\Phi_h(u),
\]
all of whose components are bounded by $N$.
Here
\begin{align*}
  |z_1-z_2-z_3+z_4|&=|\Delta(s)-\Delta(u)|+|\Delta^\tau(s)-\Delta^\tau(u)|
  \leq d\,|s-u|^\mu+\mathcal{T}^\mu d\,|s-u|^{\gamma\mu},\\
  |z_1-z_3|&=|\Delta(s)|+|\Delta^\tau(s)|\leq2\ell^\mu d,\\
  |z_1-z_2|^\delta+|z_3-z_4|^\delta
  &\leq 2\,(2\mathcal{T}^\mu R)^\delta\,|s-u|^{\delta\gamma\mu},
\end{align*}
where we used~\eqref{eq:Delta-basic}, \eqref{eq:Delta-tau}
and~\eqref{eq:Phi-holder}. Hence, since $2^{2+\delta}\leq8$,
\begin{equation}\label{eq:J2-bound}
  |J_2(s,u)|\leq M_0d\,|s-u|^\mu+M_0\mathcal{T}^\mu d\,|s-u|^{\gamma\mu}
  +8M_N\mathcal{T}^{\delta\mu}R^\delta\,\ell^\mu d\,|s-u|^{\delta\gamma\mu}.
\end{equation}

\smallskip
\noindent Then, 
by~\eqref{eq:v-decomp}, \eqref{eq:J1-bound} and~\eqref{eq:J2-bound}, $v$
satisfies the assumption of Lemma~\ref{lem:young} on every
$[u,s]\subset[a,b]$ with the exponents
$\beta,\mu,\gamma\mu,\delta\gamma\mu$, all larger than $\alpha$
by~\eqref{eq:alpha0}.
Using~\eqref{eq:young}, \eqref{eq:v-pointwise}, the definition of $c_0$ and
$s-u\leq\ell\leq1$,
\begin{align*}
  \left|\int_u^s v\,dg\right|
  &\leq c_0\Lambda(s-u)^\mu\,d\,\Bigl(2M_0\ell^\mu+2M_0\ell^{\mu+\beta}
  +M_0\ell^\mu+M_0\mathcal{T}^\mu\ell^{\gamma\mu}
  +8M_N\mathcal{T}^{\delta\mu}R^\delta\ell^{\mu+\delta\gamma\mu}\Bigr)\\
  &\leq c_0\Lambda(s-u)^\mu\,d\,\ell^{\gamma\mu}
  \Bigl(5M_0+M_0\mathcal{T}^\mu+8M_N\mathcal{T}^{\delta\mu}R^\delta\Bigr),
\end{align*}
which is~\eqref{eq:G-lip}.
\end{proof}

\subsection{A priori bound}

Consider the deterministic integral equation
\begin{equation}\label{eq:det}
\begin{aligned}
  x(t) &= \eta(0) + \int_0^t b(r,x(r),x^\tau(r))\,dr
  + \int_0^t\sigma(r,x(r),x^\tau(r))\,dg_r,
  && t\in(0,T],\\
  x(t)&=\eta(t), && t\in[-r_0,0],
\end{aligned}
\end{equation}
where $\eta\in C^{\mu}([-r_0,0];\R^d)$ is the initial path.
For $T'\in(0,T]$, a solution on $[0,T']$ is a function
$x\in C^\mu([-r_0,T'];\R^d)$ satisfying~\eqref{eq:det} for $t\in(0,T']$.
Clearly, a solution on $[0,T']$ satisfies, for
every $0\leq a<b\leq T'$,
\begin{equation}\label{eq:det-incr}
  x(t)=x(a)+F_a(x)(t)+G_a(x)(t),\qquad t\in[a,b].
\end{equation}

The following proposition gives an a priori bound for the solutions
of~\eqref{eq:det}.

\begin{proposition}\label{prop:apriori}
Define
\begin{equation}\label{eq:ell}
  \ell:=\min\Bigl\{1,\ \frac{1}{16L_0},\
  \bigl(16\,c_0M_0(1+\mathcal{T}^\mu)\Lambda\bigr)^{-1/(\gamma\mu)}\Bigr\},
  \qquad n:=\lceil T/\ell\rceil.
\end{equation}
Here we use the convention $1/0=+\infty$.
Let $T'\in(0,T]$ and let $x$ be a solution of~\eqref{eq:det} on $[0,T']$.
Then
\begin{equation}\label{eq:apriori}
  \|x\|_{\infty,[-r_0,T']}+\ell^\mu\,[x]_{\mu,[-r_0,T']}
  \leq 2^{n}\bigl(\|\eta\|_{C^\mu([-r_0,0])}+4C_1\bigr).
\end{equation}
Moreover,
\begin{equation}\label{eq:apriori-explicit}
  \|x\|_{C^{1-\alpha}([-r_0,T'])}
  \leq\bigl(1+\ell^{-\mu}\bigr)\,2^{\,1+T/\ell}\,
  \bigl(\|\eta\|_{C^{1-\alpha}([-r_0,0])}+4C_1\bigr).
\end{equation}
\end{proposition}

\begin{proof}
For $t\in[0,T']$ set
\[
  A_t:=\|x\|_{\infty,[-r_0,t]},\qquad H_t:=[x]_{\mu,[-r_0,t]},\qquad
  Z_t:=A_t+\ell^\mu H_t ,
\]
which are finite because $x\in C^\mu([-r_0,T'])$.
Let $0\leq a<b\leq T'$ with $\ell':=b-a\leq\ell$ and put $h:=[x]_{\mu,[a,b]}$.
Since $x(t)-x(a)$ is controlled by $h$ on $[a,b]$, and the H\"older
seminorm on the union of two adjacent intervals is at most the sum of the
two seminorms,
\begin{equation}\label{eq:AH}
  A_b\leq A_a+\ell^\mu h,\qquad H_b\leq H_a+h,\qquad\text{so}\qquad
  Z_b\leq Z_a+2\ell^\mu h .
\end{equation}
By~\eqref{eq:det-incr}, $h\leq[F_a(x)]_{\mu,[a,b]}+[G_a(x)]_{\mu,[a,b]}$.
We apply Proposition~\ref{prop:growth} on $[a,b]$, whose length is
$\ell'\leq\ell\leq1$, with $f=x$, so that
$\|f\|_{\infty,[-r_0,b]}=A_b$ and $[f]_{\mu,[-r_0,b]}=H_b$, and we bound
$\ell'$ by $\ell$:
\begin{equation}\label{eq:h-terms}
  h\leq \underbrace{2L_0\ell^\alpha A_b}_{\mathrm{(I)}}
  +\underbrace{\|b_0\|_{L^{1/\alpha}}}_{\mathrm{(II)}}
  +\underbrace{c_0\Lambda\bigl(K_0+M_0\bigr)}_{\mathrm{(III)}}
  +\underbrace{2c_0\Lambda M_0A_b}_{\mathrm{(IV)}}
  +\underbrace{c_0\Lambda M_0(1+\mathcal{T}^\mu)\ell^{\gamma\mu}H_b}_{\mathrm{(V)}} .
\end{equation}

The choice~\eqref{eq:ell} of $\ell$ gives the two conditions
\begin{equation}\label{eq:S12}
  2L_0\ell\leq\tfrac18,\qquad
  c_0\Lambda M_0(1+\mathcal{T}^\mu)\ell^{\gamma\mu}\leq\tfrac1{16}.
\end{equation}
The first one follows from $\ell\leq1/(16L_0)$ and the second one from the
third entry of the minimum.
Moreover, since $0<\ell\leq1$ and $\gamma\leq1$ we have
$\gamma\mu\leq\mu$ and hence
\[
  \ell^\mu\leq\ell^{\gamma\mu} .
\]
Together with $1+\mathcal{T}^\mu\geq1$, the second condition
in~\eqref{eq:S12} therefore implies
\begin{equation}\label{eq:S3}
  c_0\Lambda M_0\,\ell^\mu\leq c_0\Lambda M_0(1+\mathcal{T}^\mu)\ell^{\gamma\mu}
  \leq\tfrac1{16}.
\end{equation}
Notice that the terms (III) and (IV)
carry no power of $\ell$ in~\eqref{eq:h-terms}. Multiplying~\eqref{eq:h-terms} by $\ell^\mu$ and treating each
term separately,
\begin{align*}
  \ell^\mu\cdot\mathrm{(I)}&=2L_0\,\ell^{\mu+\alpha}A_b=2L_0\,\ell\,A_b
  \leq\tfrac18A_b
  &&\text{by $\mu+\alpha=1$ and~\eqref{eq:S12}},\\
  \ell^\mu\cdot\mathrm{(II)}&\leq\|b_0\|_{L^{1/\alpha}}
  &&\text{by $\ell\leq1$},\\
  \ell^\mu\cdot\mathrm{(III)}&=c_0\Lambda\,\ell^\mu(K_0+M_0)
  \leq\frac{K_0+M_0}{16M_0}
  &&\text{by~\eqref{eq:S3}},\\
  \ell^\mu\cdot\mathrm{(IV)}&=2\,c_0\Lambda M_0\,\ell^\mu A_b\leq\tfrac18A_b
  &&\text{by~\eqref{eq:S3}},\\
  \ell^\mu\cdot\mathrm{(V)}&\leq\tfrac1{16}\,\ell^\mu H_b
  &&\text{by~\eqref{eq:S12}}.
\end{align*}
Adding up, and recalling the definitions of $C_1$ and $Z_b=A_b+\ell^\mu H_b$,
\[
  \ell^\mu h\leq\tfrac18A_b+\tfrac18A_b+\tfrac1{16}\ell^\mu H_b+C_1
  \leq\tfrac14Z_b+C_1 .
\]
Inserting this into~\eqref{eq:AH} gives $Z_b\leq Z_a+\tfrac12Z_b+2C_1$,
that is
\begin{equation}\label{eq:one-step}
  Z_b+4C_1\leq2\,(Z_a+4C_1).
\end{equation}
Now split $[0,T']$ into at most $n$  intervals of length at
most $\ell$ and iterate~\eqref{eq:one-step}:
$Z_{T'}+4C_1\leq2^n(Z_0+4C_1)$.
Since $Z_0=\|\eta\|_{\infty}+\ell^\mu[\eta]_\mu\leq\|\eta\|_{C^\mu([-r_0,0])}$,
this is~\eqref{eq:apriori}.

Finally, from~\eqref{eq:ell},
\begin{equation}\label{eq:ell-inv}
  \frac1\ell\leq 1+16L_0+\bigl(16\,c_0M_0(1+\mathcal{T}^\mu)\Lambda\bigr)^{1/(\gamma\mu)},
  \qquad n\leq1+\frac{T}{\ell}.
\end{equation}
Since $\mu=1-\alpha$ and
$\|x\|_{\infty,[-r_0,T']}+[x]_{\mu,[-r_0,T']}
\leq(1+\ell^{-\mu})\bigl(\|x\|_{\infty,[-r_0,T']}+\ell^\mu[x]_{\mu,[-r_0,T']}\bigr)$,
the bound~\eqref{eq:apriori} and $n\leq1+T/\ell$
give~\eqref{eq:apriori-explicit}.
\end{proof}

\subsection{Existence and uniqueness: deterministic case}

We can now state and prove the existence and uniqueness result for the
deterministic equation~\eqref{eq:det}.

\begin{theorem}\label{thm:det}
Assume \textbf{(H1)} and \textbf{(H2)}, let $\alpha\in(0,\alpha_0)$ with
$\alpha_0$ as in~\eqref{eq:alpha0}, let $\tau:[0,T]\to[0,r_0]$ be
$\gamma$-H\"older, let $g:[0,T]\to\R^m$ satisfy $\Lambda_\alpha(g)<\infty$
and let $\eta\in C^{1-\alpha}([-r_0,0];\R^d)$.
Then equation~\eqref{eq:det} has a unique solution
\[
  x\in C^{1-\alpha}([-r_0,T];\R^d)\subset W_0^{\alpha,\infty}(-r_0,T;\R^d),
\]
and $x$ satisfies the bound~\eqref{eq:apriori-explicit} with $T'=T$.
\end{theorem}

\begin{proof}
By Proposition~\ref{prop:apriori}, every solution on any $[0,T']$,
$T'\leq T$, satisfies
\[
  \|x\|_{\infty,[-r_0,T']}\leq N_*,\qquad [x]_{\mu,[-r_0,T']}\leq R_*,
\]
with $N_*:=2^{n}(\|\eta\|_{C^\mu}+4C_1)$ and $R_*:=\ell^{-\mu}N_*$,
which do not depend on $T'$ or on the solution.
Define
\begin{gather*}
  R:=2\Bigl[(2L_0+2c_0\Lambda M_0)N_*+\|b_0\|_{L^{1/\alpha}}
  +c_0\Lambda(K_0+M_0)+\tfrac1{16}R_*\Bigr],\\
  N:=N_*+R,\qquad \bar R:=R_*+R,
\end{gather*}
and
\[
  \ell_*:=\min\Bigl\{\ell,\ \frac{1}{8L_N},\
  \Bigl(4c_0\Lambda\bigl(5M_0+M_0\mathcal{T}^\mu
  +8M_N\mathcal{T}^{\delta\mu}\bar R^{\,\delta}\bigr)\Bigr)^{-1/(\gamma\mu)}\Bigr\},
\]
so that
\begin{equation}\label{eq:theta}
  \Theta:=2L_N\ell_*+c_0\Lambda\,\ell_*^{\gamma\mu}
  \Bigl(5M_0+M_0\mathcal{T}^\mu+8M_N\mathcal{T}^{\delta\mu}\bar R^{\,\delta}\Bigr)
  \leq\tfrac12 .
\end{equation}

We first prove local existence.
Let $0\leq a<b\leq T$ with $b-a\leq\ell_*$ and suppose that
$x$ is a solution on $[0,a]$; if $a=0$, we take $x=\eta$.
Let
\[
  \mathcal{X}:=\bigl\{y\in C^\mu([-r_0,b];\R^d):\ y=x\text{ on }[-r_0,a],\
  [y]_{\mu,[a,b]}\leq R\bigr\},
  \qquad d(y,y'):=[y-y']_{\mu,[a,b]} .
\]
The set $\mathcal{X}$ is non-empty, since it contains the extension of $x$
by the constant $x(a)$ on $[a,b]$.
Then it is easy to check that $(\mathcal{X},d)$ is a complete metric space,
with underlying Banach space $C^\mu([-r_0,b];\R^d)$ with the norm
$\|y\|_{C^\mu([-r_0,b])}=\|y\|_{\infty,[-r_0,b]}+[y]_{\mu,[-r_0,b]}$.
Notice also that every $y\in\mathcal{X}$ satisfies, by the a priori bounds above,
\begin{equation}\label{eq:X-bounds}
  \|y\|_{\infty,[-r_0,b]}\leq N_*+R=N,\qquad
  [y]_{\mu,[-r_0,b]}\leq R_*+R=\bar R .
\end{equation}

\smallskip

Define $\mathcal{L}(y):=x$ on $[-r_0,a]$ and
$$\mathcal{L}(y)(t):=x(a)+F_a(y)(t)+G_a(y)(t)$$ for $t\in[a,b]$.

Let $y\in\mathcal{X}$ and $\ell':=b-a\leq\ell_*\leq\ell\leq1$.
We have the following bounds:
\begin{itemize}
  \item since $y=x$ on $[-r_0,a]$,
  $|x|\leq N_*$ there, and $|y(t)-y(a)|\leq R\,(t-a)^\mu$ on $[a,b]$,
  \[
    \|y\|_{\infty,[-r_0,b]}\leq N_*+\ell'^{\mu}R ;
  \]
  \item 
  $[y]_{\mu,[-r_0,b]}\leq R_*+R$, as in~\eqref{eq:X-bounds};
  \item the values of the delayed term $y^\tau$ on $[a,b]$ are values
  of $y$ on $[-r_0,b]$, so $\|y^\tau\|_{\infty,[a,b]}\leq N_*+\ell'^{\mu}R$, and
  by~\eqref{eq:Phi-holder}
  $|y^\tau(s)-y^\tau(u)|\leq\mathcal{T}^\mu(R_*+R)|s-u|^{\gamma\mu}$
  for $u,s\in[a,b]$.
\end{itemize}

Inserting the three bounds above into~\eqref{eq:F-growth}
and~\eqref{eq:G-growth} with $f=y$ gives
\begin{align*}
  [\mathcal{L}(y)]_{\mu,[a,b]}
  &\leq 2L_0\ell'^{\alpha}\bigl(N_*+\ell'^{\mu}R\bigr)+\|b_0\|_{L^{1/\alpha}}
  +c_0\Lambda(K_0+M_0)\\
  &\quad+2c_0\Lambda M_0\bigl(N_*+\ell'^{\mu}R\bigr)
  +c_0\Lambda M_0(1+\mathcal{T}^\mu)\ell'^{\gamma\mu}\bigl(R_*+R\bigr)
  \;=\;E_0+E_R ,
\end{align*}
where 
\begin{align*}
  E_0&:=2L_0\ell'^{\alpha}N_*+\|b_0\|_{L^{1/\alpha}}+c_0\Lambda(K_0+M_0)
  +2c_0\Lambda M_0N_*+c_0\Lambda M_0(1+\mathcal{T}^\mu)\ell'^{\gamma\mu}R_*,\\
  E_R&:=2L_0\,\ell'^{\alpha+\mu}R+2c_0\Lambda M_0\,\ell'^{\mu}R
  +c_0\Lambda M_0(1+\mathcal{T}^\mu)\ell'^{\gamma\mu}R .
\end{align*}
For $E_0$ we use $\ell'\leq1$ and the second condition in~\eqref{eq:S12},
which remains valid with $\ell'\leq\ell$ in place of $\ell$:
\[
  E_0\leq(2L_0+2c_0\Lambda M_0)N_*+\|b_0\|_{L^{1/\alpha}}+c_0\Lambda(K_0+M_0)
  +\tfrac1{16}R_*=\tfrac12R.
\]
On the other hand, for $E_R$ we use~\eqref{eq:S12} and~\eqref{eq:S3}, again with
$\ell'\leq\ell$ in place of $\ell$, term by term:
\begin{align*}
  2L_0\,\ell'^{\alpha+\mu}R=2L_0\,\ell'R&\leq\tfrac18R,\\
  2c_0\Lambda M_0\,\ell'^{\mu}R&\leq\tfrac18R,\\
  c_0\Lambda M_0(1+\mathcal{T}^\mu)\ell'^{\gamma\mu}R&\leq\tfrac1{16}R,
\end{align*}
so that $E_R\leq\bigl(\tfrac18+\tfrac18+\tfrac1{16}\bigr)R=\tfrac5{16}R$.
Therefore
\[
  [\mathcal{L}(y)]_{\mu,[a,b]}
  \leq E_0+E_R\leq\tfrac12R+\tfrac5{16}R=\tfrac{13}{16}R\leq R .
\]
Hence $\mathcal{L}(\mathcal{X})\subset\mathcal{X}$.
By Proposition~\ref{prop:lip}, which is applicable with the constants $N$
and $\bar R$ thanks to~\eqref{eq:X-bounds}, and by~\eqref{eq:theta},
\[
  d(\mathcal{L}(y),\mathcal{L}(y'))\leq\Theta\,d(y,y')\leq\tfrac12\,d(y,y'),
  \qquad y,y'\in\mathcal{X}.
\]
By the Banach fixed point theorem $\mathcal{L}$ has a fixed point
$y\in\mathcal{X}$. By~\eqref{eq:det-incr} for $x$ at time $a$ and the
additivity of the integrals, $y$ is a solution on $[0,b]$.

Once local existence is proved, global existence follows by an
iterative procedure.
Starting from $a=0$ and applying the local existence result successively
on the intervals
$[k\ell_*,(k+1)\ell_*\wedge T]$ we obtain, after finitely many steps, a
solution on $[0,T]$. Note that $\ell_*$ is the same at every step, because
the constants are uniform.

We finish by proving the uniqueness of the solution.
Let $x,x'$ be two solutions on $[0,T]$ and let
$t_0:=\sup\{t\in[0,T]: x=x'\text{ on }[-r_0,t]\}$.
By continuity, $x=x'$ on $[-r_0,t_0]$.
Suppose $t_0<T$. Both solutions satisfy the a priori bounds above, so by
\eqref{eq:det-incr} and Proposition~\ref{prop:lip} with $a=t_0$ and
$b=(t_0+\ell_*)\wedge T$,
\[
  [x-x']_{\mu,[a,b]}
  =\bigl[F_a(x)-F_a(x')+G_a(x)-G_a(x')\bigr]_{\mu,[a,b]}
  \leq\tfrac12\,[x-x']_{\mu,[a,b]} .
\]
Hence $x-x'$ is constant on $[a,b]$, and it vanishes at $a$; so $x=x'$ on
$[-r_0,b]$ with $b>t_0$, a contradiction. Therefore $x=x'$ on $[-r_0,T]$.
\end{proof}

\section{Proof of the main theorem}\label{sec:stoch}

Let $W^H$ be the fractional Brownian motion of Section~\ref{sec:main}.
We know that $\Lambda_\alpha(W^H)<\infty$ a.s.\ for all
$\alpha\in(1-H,\tfrac{1}{2})$, by Remark~\ref{rem:Lambda-holder}.
For $\alpha\in(1-H,\tfrac{1}{2})$ we write
$\Lambda_\alpha := \Lambda_\alpha(W^H)$.

The following integrability result is a consequence of Fernique's
theorem~\cite{Fernique1975}.

\begin{theorem}\label{thm:fernique}
For every $\alpha\in(1-H,\tfrac{1}{2})$, every $\theta\in(0,2)$ and every
$c>0$,
\begin{equation}\label{eq:fernique2}
  \E\bigl[\exp\bigl(c\,\Lambda_\alpha^{\theta}\bigr)\bigr] < \infty,
\end{equation}
and in particular $\Lambda_\alpha\in\bigcap_{p\geq 1}L^p(\Omega)$.
\end{theorem}

\begin{proof}[Proof of Theorem~\ref{thm:main}]
(a)
Let $\Omega^*$ be the set of full probability on
which~\eqref{eq:tau-holder} holds with $C_\tau(\omega)<\infty$,
$\Lambda_\alpha(W^H(\cdot,\omega))<\infty$ and
$\eta(\cdot,\omega)\in C^{1-\alpha}([-r_0,0];\R^d)$.
For every $\omega\in\Omega^*$ the assumptions of Theorem~\ref{thm:det}
are satisfied with $g=W^H(\cdot,\omega)$, $\tau=\tau(\cdot,\omega)$ and
initial path $\eta(\cdot,\omega)$.
Hence, for every $\omega\in\Omega^*$, the deterministic
equation~\eqref{eq:det} has a unique solution
\[
  X(\cdot,\omega)\in C^{1-\alpha}([-r_0,T];\R^d)
  \subset W_0^{\alpha,\infty}(-r_0,T;\R^d),
\]
and this is precisely equation~\eqref{eq:main} for that $\omega$.
This gives the existence of the solution and its uniqueness among the
processes with trajectories in $C^{1-\alpha}([-r_0,T];\R^d)$.

(b)
Let $\bar c:=\|C_\tau\|_{L^\infty(\Omega)}$, so that a.s.\
$\mathcal T(\omega)=T^{1-\gamma}+[\tau(\cdot,\omega)]_\gamma
\leq T^{1-\gamma}+\bar c=:\bar{\mathcal T}$, a deterministic constant.
Write $\mu=1-\alpha$.
By~\eqref{eq:apriori-explicit}, and using $\ell\leq1$,
$1+\ell^{-\mu}\leq 1+\ell^{-1}\leq e^{1/\ell}$ and $2^{T/\ell}\leq e^{T/\ell}$,
\[
  \|X\|_{C^{1-\alpha}([-r_0,T])}
  \leq 2\,e^{(1+T)/\ell}\bigl(\|\eta\|_{C^{1-\alpha}([-r_0,0])}+4C_1\bigr).
\]
By~\eqref{eq:ell-inv} and $\mathcal T\leq\bar{\mathcal T}$,
\[
  \frac{1}{\ell}\leq 1+16L_0+C_2\,\Lambda_\alpha^{1/(\gamma\mu)},
  \qquad C_2:=\bigl(16\,c_0M_0(1+\bar{\mathcal T}^\mu)\bigr)^{1/(\gamma\mu)},
\]
where $C_2$ is deterministic. Therefore, almost surely,
\begin{equation}\label{eq:X-bound}
  \|X\|_{C^{1-\alpha}([-r_0,T])}
  \leq C_3\,\exp\bigl(C_4\,\Lambda_\alpha^{1/(\gamma(1-\alpha))}\bigr)
  \bigl(\|\eta\|_{C^{1-\alpha}([-r_0,0])}+4C_1\bigr),
\end{equation}
with deterministic constants $C_3:=2e^{(1+T)(1+16L_0)}$ and
$C_4:=(1+T)C_2$.
The second condition in~\eqref{eq:moment-cond} says  that
\[
  \theta:=\frac{1}{\gamma(1-\alpha)}<2 .
\]
Hence, by~\eqref{eq:fernique2} and using the Cauchy--Schwarz inequality and the moment assumption on $\eta$,
\[
  \E\bigl[\|X\|_{C^{1-\alpha}}^p\bigr]
  \leq C_3^p\,\Bigl(\E\bigl[e^{2pC_4\Lambda_\alpha^{\theta}}\bigr]\Bigr)^{1/2}
  \Bigl(\E\bigl[(\|\eta\|_{C^{1-\alpha}}+4C_1)^{2p}\bigr]\Bigr)^{1/2}<\infty .
\]
The bound for $\|X\|_{\alpha,r_0}$ follows from~\eqref{eq:holder-embed}.
\end{proof}



\begin{thebibliography}{99}

\bibitem{Boukanjime2020}
B.~Boukanjime, T.~Caraballo, M.~El Fatini and M.~El Khalifi,
\textit{Dynamics of a stochastic coronavirus (COVID-19) epidemic model with
Markovian switching},
Chaos Solitons Fractals \textbf{141} (2020), 110361.

\bibitem{CarteaSanchezBetancourt2023}
\'A.~Cartea and L.~S\'anchez-Betancourt,
\textit{Optimal execution with stochastic delay},
Finance Stoch.\ \textbf{27} (2023), 1--47.

\bibitem{Fernique1975}
X.~Fernique,
\textit{Regularit\'e des trajectoires des fonctions al\'eatoires gaussiennes},
in \textit{\'Ecole d'\'Et\'e de Probabilit\'es de Saint-Flour, IV-1974},
Lecture Notes in Math.\ \textbf{480} (1975), 1--96.

\bibitem{FR2006}
M.~Ferrante and C.~Rovira,
\textit{Stochastic delay differential equations driven by fractional Brownian motion
with Hurst parameter $H > \tfrac{1}{2}$},
Bernoulli \textbf{12} (2006), 85--100.

\bibitem{FR2010}
M.~Ferrante and C.~Rovira,
\textit{Convergence of delay differential equations driven by fractional Brownian motion},
J.~Evol.\ Equ.\ \textbf{10} (2010), 761--783.

\bibitem{GarridoAtienza2011}
M.J.~Garrido-Atienza, A.~Ogrowsky and B.~Schmalfu\ss,
\textit{Random differential equations with random delays},
Stoch.\ Dyn.\ \textbf{11} (2011), 675--699.

\bibitem{HuNualart2007}
Y.~Hu and D.~Nualart,
\textit{Differential equations driven by H\"older continuous functions of
order greater than $1/2$},
in \textit{Stochastic Analysis and Applications} (Abel Symposium 2005),
Springer, 2007, 399--413.

\bibitem{KolMys}
V.~Kolmanovskii and A.~Myshkis,
\textit{Introduction to the Theory and Applications of Functional
Differential Equations},
Springer, 1999.

\bibitem{LeonTindel}
J.~Le\'on and S.~Tindel,
\textit{Malliavin calculus for fractional delay equations},
J.\ Theoret.\ Probab.\ \textbf{25} (2012), 854--889.

\bibitem{Lyons1994}
T.~Lyons,
\textit{Differential equations driven by rough signals (I):
An extension of an inequality of L.C.~Young},
Math.\ Res.\ Lett.\ \textbf{1} (1994), 451--464.

\bibitem{Mao1997}
X.~Mao,
\textit{Stochastic Differential Equations and Applications},
Horwood, 1997.

\bibitem{MaoShenYuan2008}
X.~Mao, Y.~Shen and C.~Yuan,
\textit{Almost surely asymptotic stability of neutral stochastic differential
delay equations with Markovian switching},
Stochastic Process.\ Appl.\ \textbf{118} (2008), 1385--1406.

\bibitem{Martcheva2025}
M.~Martcheva, N.~Tuncer and C.~N.~Ngonghala,
\textit{Mathematical models of epidemics with infection time},
AppliedMath \textbf{5} (2025), no.~2, 47.

\bibitem{Mohammed1984}
S.-E.A.~Mohammed,
\textit{Stochastic Functional Differential Equations},
Research Notes in Mathematics \textbf{99}, Pitman, 1984.

\bibitem{Mohammed1998}
S.-E.A.~Mohammed,
\textit{Stochastic differential systems with memory: theory, examples and applications},
in \textit{Stochastic Analysis and Related Topics VI},
Progr.\ Probab.\ \textbf{42}, Birkh\"auser, 1998, 1--77.

\bibitem{NeuenkirchNourdinTindel2008}
A.~Neuenkirch, I.~Nourdin and S.~Tindel,
\textit{Delay equations driven by rough paths},
Electron.\ J.\ Probab.\ \textbf{13} (2008), 2031--2068.

\bibitem{NR2002}
D.~Nualart and A.~R\u{a}\c{s}canu,
\textit{Differential equations driven by fractional Brownian motion},
Collect.\ Math.\ \textbf{53} (2002), 55--81.

\bibitem{Toub}
J.~Touboul,
\textit{Limits and dynamics of stochastic neuronal networks with random
heterogeneous delays},
J.\ Stat.\ Phys.\ \textbf{149} (2012), 569--597.

\bibitem{Young1936}
L.C.~Young,
\textit{An inequality of the H\"older type connected with Stieltjes integration},
Acta Math.\ \textbf{67} (1936), 251--282.

\bibitem{Zahle1998}
M.~Z\"ahle,
\textit{Integration with respect to fractal functions and stochastic calculus~I},
Prob.\ Theory Relat.\ Fields \textbf{111} (1998), 333--374.

\end{thebibliography}
\end{document}